\documentclass[10pt,reqno]{amsart}
\usepackage{epsfig,graphicx,wrapfig}
\usepackage{amssymb,amscd}

\begin{document}

\let\kappa=\varkappa
\let\eps=\varepsilon
\let\phi=\varphi
\let\mf=\mathfrak

\def\Z{\mathbb Z}
\def\R{\mathbb R}
\def\C{\mathbb C}
\def\Q{\mathbb Q}
\def\P{\mathbb P}

\def\OO{\mathcal O}
\def\CP{\C{\mathrm P}}
\def\RP{\R{\mathrm P}}
\def\conj{\overline}
\def\Beta{\mathrm{B}}
\def\O{\Omega}
\def\S{\Sigma}
\def\e{\epsilon}
\def\z{\zeta}

\renewcommand{\Im}{{\mathop{\mathrm{Im}}\nolimits}}
\renewcommand{\Re}{{\mathop{\mathrm{Re}}\nolimits}}
\newcommand{\codim}{{\mathop{\mathrm{codim}}\nolimits}}
\newcommand{\id}{{\mathop{\mathrm{id}}\nolimits}}
\newcommand{\Aut}{{\mathop{\mathrm{Aut}}\nolimits}}

\newtheorem{thm}{Theorem}[section]
\newtheorem{lem}[thm]{Lemma}
\newtheorem{prop}[thm]{Proposition}
\newtheorem{cor}[thm]{Corollary}

\theoremstyle{definition}
\newtheorem{exm}[thm]{Example}
\newtheorem{rem}[thm]{Remark}

\title{Real Analytic Sets in Complex Spaces and CR maps}
\author{Rasul Shafikov}
\date{\today}
\address{Department of Mathematics, the University of Western Ontario,
London, Canada  N6A 5B7}
\email{shafikov@uwo.ca}

\begin{abstract} If $R$ is a real analytic set in $\C^n$ (viewed as
  $\R^{2n}$), then for any point $p\in R$ there is a uniquely defined
  germ $X_p$ of the smallest complex analytic variety which contains
  $R_p$, the germ of $R$ at $p$. It is shown that if $R$ is
  irreducible of constant dimension, then the function $p\to\dim X_p$
  is constant on a dense open subset of $R$. As an application it is
  proved that a continuous map from a real analytic CR manifold $M$
  into $\C^N$ which is CR on some open subset of $M$ and whose graph
  is a real analytic set in $M\times \C^N$ is necessarily CR
  everywhere on $M$.
\end{abstract}

\maketitle

\section{Introduction}
Given a real analytic set $R$ in $\C^n$ (we may identify $\C^n$ with
$\R^{2n}$), $n>1$, we consider the germ $R_p$ of $R$ at a point $p\in
R$ and define  $X_p$ to be the germ at $p$ of the smallest (with
respect to inclusion) complex analytic set in $\C^n$ which contains
$R_p$. Then $X_p$ exists and unique for each $p\in R$, with $\C^n$
being such a set for a {\it generic} $R$. It is natural to ask how the 
dimension of $X_p$ varies with $p\in R$. Consider the following
example (Cartan's umbrella, see \cite{ca1}) 
\begin{equation}\label{ce}
  R=\left\{(z_1,z_2)=(x_1+iy_1,x_2+iy_2)\in \C^2:
  x_2(x_1^2+y_1^2)-x_1^3=0;\ y_2=0\right\}.
\end{equation}
At a point $p=(0,x_2)\in \C^2$, $x_2\ne 0$, the complex analytic
set $X_p=\{z_1=0\}$ contains $R_p$, but at the origin $X_0=\C^2$. 

The set in \eqref{ce} is irreducible but has different
dimension at different points. In particular this set is not {\it
  coherent} (see Section 2 for definitions). Our main result is that
under the assumption that $R$ is irreducible and has constant
dimension, the dimension of $X_p$ is constant on a dense open subset
of $R$. More precisely, the following result holds.   

\begin{thm}\label{t1}
  Let $\O$ be a domain in $\C^n$, $n> 1$, and $R$ be an irreducible
  real analytic subset of $\O$ of constant positive dimension. Then

  (i) There exists an integer $d>0$ and a closed nowhere dense subset
  $S$ of $R$ (possibly empty) such that for any point $p\in R\setminus
  S$, $\dim X_p =d$. 

  (ii) If $R$ is in addition coherent, then there exists a complex
  analytic set $X$, defined in a neighbourhood of $R\setminus S$ such
  that for any point $p\in R\setminus S$, the germ $X_p$ is the smallest
  complex analytic set which contains the germ $R_p$.     
\end{thm}

The function $d(p)=\dim X_p$ is upper semicontinuous on $R$, and
therefore, $\dim X_p \ge d$ for $p\in S$. We do not know any examples
where $R$ has constant dimension and $S$ is nonempty, and it would be
interesting to obtain further information about this set. The set $S$
is empty in the case when the germ of $R$ is complex analytic at some
point.    

\begin{cor}\label{tl}
  If $R\subset \O$ is an irreducible real analytic set of constant
  dimension, which is complex analytic near some point $p\in R$, then 
  $R$ is a complex analytic subset of $\O$.
\end{cor}

Again, Cartan's example $\{x_2(x_1^2+y_1^2)=x_1^3\}\subset \C^2$
provides an irreducible real analytic set, not of constant dimension,
which is complex analytic only on some part of it. Our main
application concerns CR-continuation of continuous maps whose graphs
are real analytic. 

\begin{thm}\label{t2}
  Let $M$ be a real analytic CR manifold. Let $f: M \to \C^N$ be
  a continuous map whose graph $\Gamma_f$ is a real analytic
  subset of $M\times \C^N$, $N\ge 1$. Suppose that $f$ is CR on a
  non-empty open subset of $M$. Then $f$ is a CR map. 
\end{thm}

Note that the real analytic set $\Gamma_f$ is not assumed to be
non-singular, and therefore the map $f$ need not be smooth. In this
case the condition for a continuous function to be CR is understood in
the sense of distributions. Further, even if $\Gamma_f$ is smooth, the
map $f$ may still be non-smooth, for example $f=x_2^{1/3}$ is a CR
function on $M=\{(z_1,x_2+iy_2)\in \C^2: y_2=0\}$, its graph is a
smooth real analytic set, but $f$ is not differentiable at
$\{x_2=0\}$.

\medskip

\noindent{\it Acknowledgments.} Research is partially supported by the
Natural Sciences and Engineering Research Council of Canada. The
author would like to thank S.~Nemirovski and F.~L\'arusson for helpful
discussions.  

\section{Real analytic and subanalytic sets}

In this section we briefly review basic facts about real analytic
sets. A real analytic set $R$ in an open set $\O\subset\R^n$ is
locally (i.e. in a neighbourhood of each point in $\O$) defined as the
zero locus of finitely many real analytic functions. $R$ is called
irreducible if it cannot be represented as a union of two real
analytic sets each not equal to $R$. The germ $R_p$ of a real analytic
set $R$ at a point $p\in R$ has a naturally defined complexification
$R^c_p$. If $\R^n$ is viewed as as subset of $\C^n$, then
$R^c_p\subset \C^n$ is the complex analytic germ at $p$ characterized
by the property that any holomorphic germ at $p$ which vanishes on
$R_p$ necessarily vanishes on $R^c_p$. Thus $R^c_p$ is the germ of a
complex analytic set of real dimension twice that of $R$ at $p$. To
define complexification of real analytic sets in $\C^n$ it is
convenient to introduce the following construction.    

Let $\mf{d}:\C^n_\z\to\C^{2n}_{(z,w)}$ be the map defined by
$\mf{d}(\z)=(\z,\overline \z)$. Then $\mf{D}:=\mf{d}(\C^n)$ is a
totally real embedding of $\C^n$ into $\C^{2n}$. Suppose R is a real 
analytic set of dimension $d>0$, $p\in R$, $U\subset \C^n$ is some
neighbourhood of $p$, and 
\begin{equation*}
R\cap U =\left\{\phi_j(\Re\z,\Im \z)=\phi_j(\z,\overline \z)=0, \  \
j=1,\dots,k \right\}, 
\end{equation*}
where $\phi_j(\z,\overline \z)$ are real analytic in $U$. Then the
complexification $R^c_p$ of $R_p$ can be defined as the germ at
$\mf d(p)$ in $\C^{2n}$ of the smallest complex analytic set which
contains the germ of $\mf d(R)$ at $\mf d(p)$. When $U$ and $\phi_j$
are suitably chosen, the complexification may simply be given by a 
representative  
\begin{equation*}
\left\{(z,w)\in U': \phi_j(z,w)=0\right\},
\end{equation*}
where $U'\subset \C^{2n}$ is some neighbourhood of $\mf d(p)$. Thus
$R^c_p$ is the germ of a complex analytic set of complex dimension $d$
such that $R^c_p\cap \mf{D}=\mf{d}(R_p)$. If $R_p$ is irreducible,
then so is $R^c_p$. The following proposition allows to replace germs
of analytic sets with their representatives. The proof can be found in
\cite{n}.   

\begin{prop}\label{p1}
  Let $R\subset \O$ be a real analytic set. Then for every point $p\in
  R$ there is a neighbourhood $U$ of $p$ such that  if $Q$ is a real
  analytic set in $\O$ and $R_p\subset Q_p$, then $R\cap U \subset
  Q\cap U$. 
\end{prop}

It follows from Proposition~\ref{p1} that the function $d(p)=\dim X_p$ 
is upper semicontinuous on $R$. For $R\subset\O$ real analytic denote
by $\mathcal O^{\R}(\O)$ the sheaf of germs of real analytic
functions, and by $\mathcal I(R)$ the ideal in $\mathcal O^\R(\O)$ of
germs of real analytic functions that vanish on $R$ (the so-called
sheaf of ideals of $R$). Then $R$ is called coherent if $\mathcal
I(R)$ is a coherent sheaf of $\mathcal O^\R$-modules. In fact, it
follows from Oka's theorem (which also holds in the real analytic
category) that $R$ is coherent if the sheaf ${\mathcal I}(R)$ is {\it
locally finitely generated}. The latter means that for every point
$a\in R$ there exists an open neighbourhood $U\subset \O$ and a finite
number of functions $\phi_j$, real analytic in $U$ and vanishing on
$R$, such that for any point $b\in U$, the germs of $\phi_j$ at $b$
generate the ideal $\mathcal I(R_b)$. Note that the corresponding
statement for complex analytic sets always holds by Cartan's theorem,
i.e. every complex analytic set is coherent.  

\begin{prop}\label{p2}
  Let $\O\subset\C^n$ be an open set, and let $R\subset \O$ be an
  irreducible real analytic set of constant positive
  dimension $d$. Suppose  $R^c_a=A_{\mf d(a)}$, where $A_{\mf d(a)}$
  is a germ at point $\mf d(a)$ of some irreducible complex analytic
  set $A$ defined in some open set in $\C^{2n}$. Then for any point
  $b\in R$ sufficiently close to $a$, 
  
  (i) $R^c_b\subset A_{\mf d(b)}$, where $A_{\mf d(b)}$ is the germ of
  $A$ at $\mf d(b)$. Further, $\dim R^c_b=\dim A_{\mf d(b)}$, and thus
  $R^c_b$ is the union of certain irreducible components of $A_{\mf
  d(b)}$. 

  (ii) If $R$ is coherent, then $R^c_b=A_{\mf d(b)}$.
\end{prop}

\begin{proof} (i) Let $A$ be an irreducible complex analytic subset of
  some open set $U'\subset \C^{2n}$, $\mf d(a)\in U'$, such that
  $A_{\mf d(a)}=R^c_a$. It follows from Proposition~\ref{p1} that if
  $U'$ is sufficiently small, then $\mf d(R)\cap U' \subset
  A$. In particular, if $b\in R$ is close to $a$, then $\mf
  d(R_b)\subset A$, and therefore, $R^c_b\subset A$. Now, since  
  $A$ and $R^c_b$ both have dimension $d$ and contain $\mf d(R_b)$,
  which is a totally real subset of real dimension $d$, it follows
  that $R^c_b$ must coincide with the union of some irreducible
  components of $A_{\mf d(b)}$.

(ii) Since $R$ is coherent, there exist functions
  $\phi_1(\z,\overline \z),\dots, \phi_k(\z,\overline \z)$ real
  analytic in some open set $U\subset \C^n$ containing $a$, such that
  the germs of these functions at any point $b\in U$ generate the ideal
  ${\mathcal I}(R_b)$. Then $R^c_b$ is defined by the equations
  $\phi_j(z,w)=0$. By the uniqueness theorem for complex analytic
  sets, $R^c_b=A_b$. 
\end{proof}

In Proposition~\ref{p2}(ii) the assumption that $R$ is coherent cannot
be in general replaced by the assumption that $R$ has constant
dimension. Indeed, consider the set (cf. \cite{ca1}) 
\begin{equation}\label{c2}
R=\{x\in\R^3: x_3(x_1^2+x_2^2)(x_1+x_2)=x_1^4\}.
\end{equation}
This is an irreducible real analytic set of constant dimension which
is not coherent at the origin. We may naturally identify $\R^3$ with
the set $\C^3\cap\{y_1=y_2=y_3=0\}$. Then the set
\begin{equation*} 
X=\{z=x+iy\in \C^3: z_3(z_1^2+z_2^2)(z_1+z_2)=z_1^4\}
\end{equation*}
can be viewed as the complexification of $R$ at the origin. However,
at any point $p=(0,0,x_3),\ x_3\ne 0,$ the germ $X_p$ is reducible,
and only one of its components is the complexification of $R_p$, which
is irreducible. The set $R$ in \eqref{c2} (viewed as a subset of
$\C^3$) is also an example of an irreducible non-coherent real
analytic set which has constant dimension and such that there is no
globally defined $X$ such that $X_p$ is the smallest complex analytic
germ containing $R_p$ for all $p\in R$. 

In the proof of Theorem~\ref{t2} we will use some results concerning
subanalytic sets. A subset $R$ of a real analytic manifold $M$ is
called {\it semianalytic} if for any point $p\in M$ there exist a
neighbourhood $U$ and a finite number of functions $\phi_{jk}$ and
$\psi_{jk}$ real analytic in U such that  
\begin{equation*}
  R\cap U =\cup_{j} \{\z\in U: \phi_{jk}(\z)=0,\ \psi_{jk}(\z)>0,\
  k=1,\dots,l\}. 
\end{equation*}
In particular, a real analytic set is semianalytic. A subset $S$ of a
real analytic manifold $M$ is called {\it subanalytic} if for any
point in $M$ there exists a neighbourhood $U$ such that $S\cap U$ is a  
projection of a relatively compact semianalytic set, that is there
exists a real analytic manifold $N$ and a relatively compact
semianalytic subset $R$ of $M\times N$ such that $S\cap U=\pi(R)$,
where $\pi:M\times N\to M$ is the projection.

The main reason for introducing the class of subanalytic sets comes
from the fact that the images of semianalytic (in particular real
analytic) sets under real analytic maps are subanalytic. Semi- and
subanalytic sets enjoy many properties of real analytic sets, for
example, a finite union, intersection and set-theoretic complement of
such sets is again in the same class. Further, semi- and subanalytic
sets admit stratifications satisfying certain properties.  

Given a subanalytic set $S$, we call a point $p\in S$ {\it regular} if
near $p$ the set $S$ is just a real analytic manifold of dimension
equal to that of $S$ (i.e. maximal possible). Denote by ${\rm reg}(S)$
the set of all regular points. Points which are not regular form the
set of singular points, ${\rm sing}(S)$. We will use the following
result due to Tamm, \cite{ta}, Theorem 1.2.2(v).   

\begin{prop}[\cite{ta}]\label{t} If $S$ is a subanalytic set, then
  ${\rm reg}(S)$ and ${\rm sing}(S)$ of $S$ are both
  subanalytic. Moreover, $\dim {\rm sing}(S)\le \dim S-1$, unless
  $S=\emptyset$.  
\end{prop}

For more about semianalytic and subanalytic sets see e.g. \cite{bm} or
\cite{hi}. A detailed discussion of real analytic sets can be found in 
\cite{ca1}, \cite{m} or \cite{n}.

\section{Pre-images of projections of analytic sets}

Let $R\subset\C^n$ be a real analytic set, $0\in R$, and $R$
irreducible at the origin. Let $A\subset \C^{2n}$ be a representative
of $R^c_0$, the complexification of the germ $R_0$, and let $\mf
d:\C^n_\zeta\to \C^{2n}_{(z,w)}$ be the totally real embedding. After
rescaling we may assume that $A$ is an irreducible complex analytic
subset of the unit polydisc 
$\Delta^{2n}=\{(z,w)\in\C^n\times\C^n : |z_j|<1,\ |w_j|<1,\
j=1,\dots,n\}$, and $\mf d(R)\cap \Delta^{2n}\subset A$. Let 
\begin{equation}\label{proj}
\pi : \Delta^{2n}_{(z,w)}\to \Delta^n_{z}
\end{equation}
be the coordinate projection onto $z$-subspace. For the proof of
Theorem~\ref{t1} we will need the following result. 

\begin{lem}\label{l} 
  There exist a closed nowhere dense subset $S$ of $A$, a
  neighbourhood $U$ of $A\setminus S$ in $\Delta^{2n}$, and a complex
  analytic subset $Y$ of $U$ with the following properties:   
  
  (a) $S$ does not divide $A$, and $\mf d(R)\cap S$ is nowhere
  dense in $\mf d(R)$.

  (b) $Y$ can be locally given as the zero locus of a system of
  holomorphic equations each of which is independent of the variable
  $w$. 

  (c) $A\cap U\subset Y\subset \pi^{-1}(\pi(A))$.
\end{lem}

\begin{proof} For $1\le j\le n$ define
\begin{equation}
\pi_j:\C^{2n}_{(z,w)}\to \C^{2n-1}_{(z,w_1,\dots, w_{j-1},w_{j+1},
  \dots,w_n)} 
\end{equation}
to be the coordinate projection parallel to $w_j$-direction, and
observe that for any $1\le k\le n$,
\begin{equation}\label{proj-incl}
\pi^{-1}_k\circ \pi_k\circ \pi^{-1}_{k-1}\circ \pi_{k-1}\dots
\pi^{-1}_1 \circ\pi_1 (A) \subset \pi^{-1}(\pi(A)).
\end{equation}
For $p\in A$ let $l_p\pi_1$ denote the germ of the fibre of $\pi_1|_A$
at $p$, i.e. the germ at $p$ of the set $\pi^{-1}_1(\pi_1(p))\cap
A$. Then by the Cartan-Remmert theorem (see e.g. \cite{l}) the set  
\begin{equation}
E_1=\{p\in A: \dim l_p\pi_1>0\}
\end{equation}
is complex analytic. Suppose that $\dim E_1 = \dim A$. Then, since $A$
is irreducible, $E_1=A$, and therefore every point of $A$ has a fibre
of dimension one. It follows then that $\pi^{-1}_1(\pi_1(A))=A$, and
so $\pi^{-1}_1(\pi_1(A))$ is a complex analytic subset of
$\Delta^{2n}$.   

Assume now that $\dim E_1 < \dim A$. Then $\mf d(R)\not\subset E_1$,
since otherwise $\mf d(R)$ would be contained in a complex analytic
set of dimension smaller than that of $A$, but $A$ is a representative
of the complexification of $R_0$, which means that $A$ is the smallest
complex analytic set containing $\mf d(R)$. Since every point in $E_1$
has a fibre of dimension one, $\pi^{-1}_1(\pi_1(E_1))=E_1$, and if
$p\in A\setminus E_1$, then $\pi_1(p)\notin \pi_1(E_1)$. 

Suppose that $p_0\in A\setminus E_1$. Then $\pi^{-1}_1(\pi_1(p_0))\cap
A$ is a finite set. Therefore, there exists a neighbourhood
$\Omega_0\subset \Delta^{2n}$ of $p_0$ such that the restriction of
$\pi_1$ to $\Omega_0\cap A$ is a finite map. Furthermore, we can
choose $\O_0$ to be of the form $\Omega_0=\O'_0\times \O''_0\subset
\Delta^{2n-1}\times \Delta_{w_1}$, with the property that $A\cap(\O'_0
\times (\partial\O''_0)) =\varnothing$. It follows that
$\pi_1|_{A\cap\O_0}$ is a proper map, and by the Remmert proper
mapping theorem, $\pi_1(A\cap \O_0)$ is a complex analytic subset of
$\Omega'_0$. Hence, $\pi^{-1}_1(\pi_1(A\cap \O_0))$ is a complex
analytic subset of $\Omega_0$. Furthermore, 
\begin{equation}\label{y1}
\pi^{-1}_1\left(\pi_1 (A\cap \Omega_0 )\right) \cap \Omega_0= 
\left((\pi_1(A\cap\Omega_0))\times \Delta_{w_1}\right) \cap \Omega_0
\end{equation}
defines a complex analytic set that can be given by a system of
equations independent of $w_1$. We set $S_1=E_1$, and thus we proved 
that there exists a neighbourhood $U_1$ of $A\setminus S_1$ and a
complex analytic subset $Y_1$ of $U_1$ that locally can be
represented as in \eqref{y1}. It follows from \eqref{proj-incl} that
$Y_1\subset \pi^{-1}(\pi(A))$.  

We now argue by induction. Suppose that there exist a closed nowhere
dense subset $S_k$ of $A$, which satisfies condition (a) of the lemma,
and a complex analytic subset $Y_k$ of some neighbourhood $U_k$ of
$A\setminus S_k$, which satisfies (c) and can be locally given in the
form $\{\phi(z, w_{k+1},\dots w_n)=0\}$, $k<n$. If the set $Y_k$ is
reducible, we keep only one irreducible component of $Y_k$ which
contains $A\cap U_k$, and for simplicity denote this component again
by $Y_k$. (Observe that $A\cap U_k$ is still irreducible, since the
regular part of $A\cap U_k$ is connected.) Consider the set   
\begin{equation}
E_{k+1}=\{p\in Y_k: \dim l_p\pi_{k+1}>0\}.
\end{equation}
Then as before,
$E_{k+1}$ is a complex analytic subset of $Y_k$. If $\dim E_{k+1} =
\dim Y_k$, then  $E_{k+1}=Y_k=\pi^{-1}_{k+1}(\pi_{k+1}(Y_k))$, and
\begin{equation*}
\pi^{-1}_{k+1}(\pi_{k+1}(Y_k))=\left( \pi^{-1}_{k+1}(\pi_{k+1}(Y_k))
\cap \{(z,w):w_{k+1}=0\}\right) \times \Delta_{w_{k+1}}.
\end{equation*}
This show that $Y_k$ is a complex analytic subset of $U_k$ that
satisfies (c), and can be given by a system of equations independent
of $(w_1,\dots, w_{k+1})$.  

Suppose now that $\dim E_{k+1} < \dim Y_k$. If $A\cap U_k\subset
E_{k+1}$, then we simply define $Y_{k+1}=E_{k+1}$. This defines a
complex analytic subset of $U_k$ with all the required properties. 
The remaining case is $A\not\subset E_{k+1}$. Then every point $a\in A
\setminus E_{k+1}$ is contained in some neighbourhood in which the map 
$\pi_{k+1}|_{Y_{k}}$ is proper, and we may repeat the argument above
to define $Y_{k+1}$ to be a complex analytic subset of some
neighbourhood of $A\setminus S_{k+1}$, where $S_{k+1}=S_k\cup (A\cap
E_{k+1})$, such that locally $Y_{k+1}$ is given as
$\pi^{-1}_{k+1}(\pi_{k+1}(Y_k))$.
By assumption, $S_{k+1}$ is closed nowhere dense in $A$,
does not divide $A$, and does not contain $\mf d(R)$. By construction,
$Y_{k+1}$ is locally given as $\pi_{k+1}(Y_k)\times \Delta_{w_{k+1}}$,
and therefore it can be defined by a system independent of $(w_1,\dots,
w_{k+1})$. Finally, by \eqref{proj-incl}, $Y_{k+1}\subset
\pi^{-1}(\pi(A))$. 

After $n$ steps the set $Y_{n}$ defined in a neighbourhood of
$A\setminus S_n$ will satisfy the lemma. 
\end{proof}  

\section{Proof of Theorem \ref{t1}}

For each $p\in R$ denote by $d(p)$ the dimension of the germ
$X_p$. Since $d(p)$ is a positive integer-valued function on 
$R$, there exists a minimum, say $d$. Let $q$ be a point where $d(p)$
attains its minimum, and let $X_q$ be the germ of a complex analytic
set at $q$ with $\dim X_q = d$ which contains $R_q$. Furthermore,
since the set of points where $R$ is not locally irreducible is
contained in the singular part of $R$ (and hence nowhere dense in
$R$), the point $q$ can be chosen in such a way that $R$ is locally
irreducible at $q$. This implies that $X_q$ can also be chosen to be
irreducible. Let $\Omega$ be a connected open neighbourhood of $q$, and
$X\subset \Omega$ be a particular representative of the germ $X_q$ such
that $X_p$ is the germ of the smallest dimension containing $R_p$ for
all $p\in R\cap \Omega$, but $X$ is either not defined at some point
$p_0 \in \partial \Omega \cap R$ or $X_{p_0}$ does not satisfy the
described property. We consider two cases depending on whether
$R_{p_0}$ is irreducible or not.   

Along with the complexification of $R$ we may also consider the
complexification in $\C^{2n}$ of the set $X$. If $\{\phi_k(\z)=0\}$ is
the system of holomorphic equations defining $X$, then the system
$\{\phi_k(z)=0\}$ (variables $w_j$ are not involved) defines a complex
analytic set $X_z$ on some open set in $\C^{2n}$. We note that
$X_z\cap\mf{D} = \mf{d}(X)$, and that the canonical complexification
of $X$ (viewed as a real analytic set) can be recovered from $X_z$ as
\begin{equation}
X^c=X_z\cap X_w :=\{\phi_k(z)=0\}\cap \{\overline \phi_k(w)=0\}.
\end{equation}
In particular, $X^c \subset X_z$ on a non-empty open set in $\C^{2n}$ 
where $X^c$ and $X_z$ are both well-defined. 

Conversely, if $Y$ is a complex analytic set in $\C^{2n}$ defined by a
system of equations $\phi_k(z)=0$, which are independent of $w$, then
$Y$ induces a complex analytic set $X$ in $\C^n$ that can be defined
as  
\begin{equation*}
X:=\mf d^{-1}\left(\mf D \cap \{\phi_k(z)=0\} \cap \{\overline \phi_k
(w)=0\}\right). 
\end{equation*}
If fact, 
$\mf D \cap \{\phi(z)=0\} \cap \{\overline \phi (w)=0\} = \mf D \cap
Y$, and $X=\mf d^{-1}(\mf D\cap Y)$.

Consider first the case when $R_{p_0}$ is irreducible. After a
translation we may assume that $p_0=0$, and so any neighbourhood of
the origin contains an open piece of $X$. To simplify the exposition
we denote by $A$ some representative of $R^c_0$ and assume without
loss of generality that $A$ is a complex analytic subset of the
polydisc $\Delta^{2n}$ defined by a system of equations holomorphic in 
$\Delta^{2n}$. Here $(z,w)=(z_1,\dots,z_n, w_1,\dots,w_n)\in \mathbb
C^n \times \mathbb C^n$, and the complexification of $R$ comes from
identifying $\z$ with $z$, and $\overline \z$ with $w$. The set $A$
can be chosen irreducible. Let $\pi$ be as in \eqref{proj}. By
Lemma~\ref{l} there exist a closed nowhere dense subset $S\subset A$,
which does not divide $A$ and $\mf d(R)\not\subset S$, a neighbourhood
$U$ of $A\setminus S$ in $\Delta^{2n}$, and a complex analytic subset
$Y$ of $U$, which may be locally defined by a system of equations
independent of $w$, such that $A\cap U\subset Y \subset
\pi^{-1}(\pi(A))$. Since $U$ is connected, and $Y$ can be chosen
irreducible, we may assume that $Y$ has constant dimension.  

Let $p\in R\cap X$ be arbitrarily close to the origin, and $\mf
d(p)\in \Delta^{2n}\setminus S$. By Proposition~\ref{p2}(i), there
exists a neighbourhood $V\subset \Delta^{2n}$ of $\mf d(p)$ such that
certain components of $A\cap V$ coincide with $R^c_p\cap V$. Denote
them by $\tilde A$. Since $R^c_p\subset X^c \subset X_z$, we conclude
that $\tilde A\subset X_z$. Therefore, $\pi^{-1}(\pi(\tilde A))\subset
X_z$. Indeed,  
\begin{equation*} 
\pi^{-1}(\pi(\tilde A))\subset \pi^{-1}(\pi(X_z))=X_z, 
\end{equation*}
where the last equality holds because $X_z$ is defined by a system of
equation independent of $w$. On the other hand, $X$ is the smallest
complex analytic germ containing $R$, and therefore,
$\pi^{-1}(\pi(\tilde A))= X_z$, as otherwise, the set
$\pi^{-1}(\pi(\tilde A))$ would induce a smaller complex analytic set
in $\C^n$ containing $R_p$. 

We now claim that $\dim Y = \dim X_z$. First, observe that
\begin{equation}\label{XA}
\dim \pi(R^c_p) = \dim \pi(A).
\end{equation}
Indeed, suppose that on the contrary, $\dim \pi(R^c_p)<\dim
\pi(A)$. Let $k$ be the generic dimension of the fibre
$\pi^{-1}(z)\cap A$ for $z\in\pi(A)$. Then since $\dim R^c_p=\dim A$,
\begin{equation*}
  R^c_p \subset \{a\in A : \dim l_a \pi > k\}.
\end{equation*}
The latter is a complex analytic subset of $A$ by the Cartan-Remmert
theorem, and by the assumption, it is a proper subset of $A$. But this 
contradicts irreducibility of $A$. Thus \eqref{XA} holds, which
implies $\dim \pi^{-1}(\pi(R^c_p))= \dim Y$. But $\dim
\pi^{-1}(\pi(R^c_p))=\dim X_z$, as otherwise, $\pi^{-1}(\pi(R^c_p))$
induces a complex analytic set in $\C^n$ which contains $R_p$ and
which is smaller than the set induced by $X_z$. This proves the
claim. Finally, since $\pi^{-1}(\pi(\tilde A))=X_z$, the set $Y\cap V$
contains $X_z\cap V$ as a union of locally irreducible components at
$p$.    

Thus we proved that if $R$ is locally irreducible at $p_0=0$, then 
for any point $p\in R\setminus \mf d^{-1}(\mf d(R)\cap S)$, the set
$\mf d^{-1}(Y\cap \mf D)$ defines a complex analytic germ in $\C^n$ at 
$p$ which has dimension $d$ and contains $R_{p}$. To complete the
proof of part (i) of Theorem~\ref{t1} it remains to consider the case
when $R_{p_0}$ is reducible. Note that the above construction produces
a complex analytic germ of dimension $d$ which contains a dense open
subset of one of the irreducible components of $R_{p_0}$. 

We claim that given any two points on $R$ there exists a path
$\gamma\subset R$ that connects these points and satisfies the
property that if $a\in \gamma$ is a point where $R$ is locally 
reducible, then $\gamma$ stays in the same local irreducible component
of $R$ at $a$. Arguing by contradiction, denote by $\Sigma$ the set of
all points on $R$ that can be connected with a given point $q$ by a
path $\gamma$ which satisfies the above property, and suppose that
$\Sigma\ne R$. We claim that $\Sigma$ is a real analytic set. Indeed,
if $p\in \Sigma$ is a smooth point of $R$, then clearly a full
neighbourhood of $p$ in $R$ is contained in $\Sigma$. If $p\in \Sigma$
is a singular point of $R$, then either $R_p$ irreducible, in which
case again a full neighbourhood of $p$ in $R$ is contained in
$\Sigma$, or $R_p$ is reducible, and only some components of $R_p$
are in $\Sigma$. In any case, $\Sigma$ is a real analytic subset of
some neighbourhood of $p$. Since $\Sigma$ is clearly closed, it
follows that it is a real analytic set. Let $\Sigma':=
\overline{R\setminus \Sigma}$, where the closure is taken in $R$. Then
$\Sigma'$ is also a real analytic set. Indeed, if $a\in R\setminus
\Sigma$, then a full neighbourhood of $a$ in $R$ is in $\Sigma'$, and
if $a\in (\overline{R\setminus \Sigma})\setminus (R\setminus \Sigma)$,
then $a\in R$ is a point at which $R_a$ is reducible, and therefore
near $a$ the set $\Sigma'$ coincides with some irreducible components
of $R_a$. Since $\Sigma'$ is analytic near any of its points and
closed, it follows that $\Sigma'$ is a real analytic set. Thus
$R=\Sigma\cup \Sigma'$, but this contradicts irreducibility of
$R$. Hence, $\Sigma=R$, and that proves the claim. 

So if in the situation above $R_{p_0}$ is reducible, we find a path 
$\gamma$ with the described property which connects $q$ with the
points on other components of $R_{p_0}$ and repeat the above
construction along $\gamma$ sufficiently many times. This proves that
there exists a dense open subset $U$ of some neighbourhood of $p_0$
such that for any $p\in U$, the germ $R_{p}$ is contained in some
complex variety of dimension $d$.   

For the proof of part (ii) simply observe, that by
Proposition~\ref{p2}(ii), $A=X^c$ near $p$, and therefore, $Y$ defines  
analytic continuation of the set $X_z$ to a neighbourhood of the
origin. This in its turn provides analytic continuation of the set $X$ 
to a neighbourhood of a dense open subset of $R$. Further, the same
holds if we repeat the above construction near any other point on $R$,
and therefore there exists a complex analytic set $X$ in a
neighbourhood of a dense open subset of $R$ with the desired
properties. Theorem~\ref{t1} is proved.  

\begin{proof}[Proof of Corollary \ref{tl}]
It follows from Theorem~\ref{t1} that there exists an open set
$\Omega\subset \C^n$ such that $R\cap \Omega$ is a complex analytic
set which is dense in $R$. We now use the result of Diederich and
Forn\ae ss \cite{df}. The claim in the proof of Theorem~4 in~\cite{df}
states that if $X_q$ is a complex analytic germ at $q\in R$, and
$X_q\subset R$, then there exists a neighbourhood $U$ of $q$,
independent of $X_q$, such that $X_q$ extends to a closed complex
analytic subset of $U$. It follows that $\Omega$ can be chosen to be a 
neighbourhood of $R$, which proves that $R$ is a complex analytic
set. 
\end{proof}

\section{Proof of Theorem \ref{t2}}

First note that $\Gamma_f$ is an irreducible real analytic set of
constant dimension. Let $\tilde M$ be the subset of $M$ on which $f$
is CR. Since the set of singular points of the real analytic set
$\Gamma_f$ is nowhere dense in $\Gamma_f$, there exists a point $p\in
\tilde M$ such that $(p,f(p))$ is a smooth point of
$\Gamma_f$. Further, the point $p$ can be chosen in such a way that
$\Gamma_f$ near $(p,f(p))$ is the graph of a smooth map on $M$. In
fact, by the real analytic version of the implicit function theorem
(see e.g. \cite{ca}) there exists a neighbourhood $U_p$ of $p$ such
that the map $f: M\cap U_p \to \C^N$ is real-analytic.  

If the CR codimension of $M$ is zero, then $M$ is simply a complex
manifold, and the map $f$ is holomorphic in $U_p$. Therefore,
$\Gamma_f$ is complex analytic over $U_p$, and by Corollary~\ref{tl}
is a complex analytic set in $M\times \C^N$. Since the projection
$\pi: \Gamma_f\to M$ is injective, it follows (see e.g. \cite{c}) that
$\Gamma_f$ is itself a complex manifold, $\pi: \Gamma_f\to M$ is
biholomorphic, and therefore, $f=\pi'\circ \pi^{-1}$ is holomorphic
everywhere on $M$ (here $\pi':\Gamma_f\to \C^N$ is another
projection). If the CR dimension of $M$ is zero, then there is nothing to
prove since any function is CR. Hence, we may assume that both the CR
dimension and the CR codimension of $M$ are positive.   

The problem is local, therefore, it is enough to prove that $f$ is CR
in a neighbourhood of a point $q\in M$ which is a boundary point of
$\tilde M$, and then use a continuation argument. By \cite{ah}, there
exists a neighbourhood $U$ of $q$ in $M$ such that $U$ can be
generically embedded into $\C^n$ for some $n>1$, i.e. $n$ is the sum
of the CR-dimension and codimension of $M$. Thus, without loss of
generality we may assume that $M$ is a generic real analytic
submanifold of some domain in $\C^n$, and $f:M\to\C^N$ is a continuous
map which is a real analytic CR map on some non-empty subset $\tilde
U$ of $M$. By \cite{t}, every component $f_j$ of $f|_{\tilde U}$
extends to a function $F_j$ holomorphic in some neighbourhood of
$\tilde U$. Then the map $F=(F_1,\dots,F_N)$ defines a complex
analytic set in $\C^n\times\C^N$ of dimension~$n$, namely, its graph
$\Gamma_F$. By construction $\Gamma_F$ contains the set $\Gamma_f$.   

Observe that $\Gamma_F$ is the smallest complex analytic set which
contains $\Gamma_f$. Indeed, suppose, on the contrary, that  there
exists a complex analytic set $A$, $\dim A < n$, which contains a
non-empty subset of $\Gamma_f$. Let $\pi: \C^{n+N}\to \C^n$ be the
projection. Then $M\subset\pi(A)$, where $\pi(A)$ is a countable union   
of complex analytic sets in $\C^n$ of dimension at most $n-1$. This is
however impossible, since $M$ is a generic submanifold of $\C^n$. 

We now show that $f$ is CR everywhere on $M$. By Theorem \ref{t1},
there exist a closed nowhere dense set $S\subset \Gamma_f$, and a 
complex analytic subset $X$ of a neighbourhood of $\Gamma_f\setminus
S$, $\dim X=n$, which contains $\Gamma_f\setminus S$. Suppose first
that $p\in M$, and $(p,p')\in \Gamma_f\setminus S \subset \C^n\times
\C^N$. Choose neighbourhoods $\O$ and $\O'$ of $p$ and $p'$
respectively, so that $X\cap (\O\times \O')$ is complex analytic. Let
$\pi: X\to \O$ be the projection, and let  
\begin{equation*}
  E=\{x\in X: \dim l_x \pi >0\}.
\end{equation*}
Let $k$ be the  generic dimension of $\pi^{-1}(\z)\cap X$ for
$\z\in\pi(X)$. Then $k=0$, since otherwise, $\pi(X)$
is a locally countable union of complex analytic sets of dimension at
most $n-k<n$, and therefore $\pi(X)$ cannot contain a generic
submanifold $M$. Therefore, $\dim E < n$, and in particular,
$M\not\subset \pi(E)$. Consider the set $E\cap\Gamma_f$. This is a
real analytic subset of $\O\times \O'$, and therefore, its projection, 
$T:=\pi(E\cap\Gamma_f)$, is a subanalytic set in $\O$. From the above
considerations, $T\ne M$.  

We first show that for any point $\z\in (M\cap \O)\setminus T$, the
map $f$ is CR at $\z$. Since $\pi(E)$ is closed, there exists a
neighbourhood $V\subset\O$ of $\z$ such that $V\cap
\pi(E)=\varnothing$. We may further shrink $V$ and choose a   
neighbourhood $V'$ of $f(\z)$ such that the projection $\pi:
X\cap(V\times V')\to V$ is proper. In particular this implies that
$\pi$ is a branched covering. Let $G\subset X\cap (V\times V')$ be the  
branch locus of $\pi$. Then near any point $a\in (M\cap V)\setminus
\pi(G)$ the map $\pi^{-1}: V\to X$ splits into a finite number of
holomorphic maps. It follows that there exists a branch of $\pi^{-1}$,
say $\sigma$, such that near $a$, the map $f$ coincides with the
restriction to $M$ of a holomorphic map $\pi'\circ \sigma:V \to
\C^N$. Therefore $f$ is CR. The set $M\cap\pi(G)$ is the intersection
of a real analytic manifold and a complex analytic subset of $V$, and
therefore it admits a stratification into a finite number of smooth
components. Each of these components is a removable CR-singularity for
$f$. More precisely one has the following result. 

\begin{lem}\label{dpl}
  Let $M$ be a smooth generic submanifold of $\C^n$, of positive CR
  dimension and codimension. Let $S\subset M$ be a smooth submanifold
  with $\dim S < \dim M$. Then any function $h$ continuous on $M$ and
  CR on $M\setminus S$ is CR on $M$.
\end{lem}

This is a trivial generalization of Proposition 4 in \cite{dp} where
the result is stated for real hypersurfaces. The proof is the same,
with the only difference being the degree of the form $\phi$. Applying
Lemma~\ref{dpl} to each smooth component of $M\cap\pi(G)$ we deduce
that $f$ is CR on $M\setminus T$. 

To prove that $T$ is a removable CR-singularity for $f$ we observe
that by Proposition~\ref{t}, the regular part of $T$, is a smooth
submanifold of $M$, and therefore, by Lemma~\ref{dpl}, ${\rm reg}(T)$
is a removable singularity for the map $f$. The set ${\rm sing}(T)$ is  
subanalytic of dimension strictly less than that of $T$, and we may
repeat the process by induction. After finitely many steps, we
conclude that $f$ is CR near $p$. 

To complete the proof of the theorem it remains to consider the case 
$(p,p')\in \Gamma_f\cap S$. We recall the construction in
Lemma~\ref{l}. Let $A$ be some representative of the complexification
of the germ $(\Gamma_f)_{(p,p')}\subset \C^{2n+2N}$. In the
notation of Lemma~\ref{l}, let $k$ be the smallest integer such that
(i) $Y_k$ defines a complex analytic subset in a neighbourhood of 
$A\setminus \tilde S$, where $\tilde S = S_1\cup \dots \cup
S_k$, and $S_j$ are defined as in the proof of Lemma~\ref{l}, (ii)
$A\subset Y_k\subset \pi^{-1}_{(w,w')}(\pi_{(w,w')}(A))$, where
$\pi_{(w,w')}:\C^{2n+2N}_{(z,z',w,w')}\to \C^{n+N}_{(z,z')}$, and
(iii) $Y_k$ is locally defined by a system of equations independent of
$(w,w')$. Then for $\mf d: \C^{n+N} \to \C^{2n+2N}$,
$\mf d(S)\subset\tilde S$ near $\mf d(p,p')$, and for any
point $q\in A\setminus \tilde S$, there exists a small neighbourhood
of $q$ where $Y_k=\{\phi_\nu(z,z')=0\}$, for some $\phi_\nu(z)$
holomorphic near $q$. Let    
\begin{equation}\label{a}
a\in S_k\setminus (\bigcup_{j=1}^{k-1}S_j)\cap \mf d(\Gamma_f).
\end{equation}
Then there exists a neighbourhood $U$ of $\mf d^{-1}(a)$ such that 
near any point $b\in \mf d^{-1}((A\setminus S_k)\cap \mf D)\cap U
\subset \Gamma_f$, the set $\Gamma_f$ is contained in a complex
analytic set of dimension $n$, and $\mf d^{-1}(S_k \cap \mf
d(\Gamma_f))$ is a real analytic subset of $\Gamma_f\cap U$. We now may
repeat the argument which we used to prove that $T$ is a removable
CR-singularity for $f$. Indeed, it follows that $f$ is CR in
$\pi(((\Gamma_f \setminus \mf d^{-1}(S_k)))\cap U)$,  and
$\pi(\Gamma_f\cap \mf d^{-1}(S_k))$ is a subanalytic set in $M$. Using
Proposition~\ref{t} and Lemma~\ref{dpl} we show that $f$ is CR on
every smooth component of $\pi(\Gamma_f\cap \mf d^{-1}(S_k))$. 

This shows that $\Gamma_f$ is the graph of a CR map at every point of
$\mf d^{-1}(S_k)\cap \Gamma_f$. By construction, the set $S_{k-1}$ is
complex analytic, and therefore, $\mf d^{-1}(S_{k-1})\cap \Gamma_f$ is
real analytic, and therefore the same procedure as before
applies. Arguing by induction we show that for all $j=1,...,k$ the set
$\pi(\Gamma_f\cap \mf d^{-1}(S_j))$ is a removable CR-singularity for
$f$. This completes the proof of Theorem~\ref{t2}.  



\bigskip
{\small

\begin{thebibliography}{99}

\bibitem{ah} A. Andreotti and C. Denson Hill, {\it Complex
  characteristic coordinates and tangential Cauchy-Riemann equations.}
  Ann. Scuola Norm. Sup. Pisa (3)  {\bf 26}  (1972), 299--324.

\bibitem{bm} E. Bierstone and P. Milman, {\it Semianalytic and
  subanalytic sets.} Inst. Hautes Études Sci. Publ. Math. No. 67,
  (1988), 5--42. 

\bibitem{ca1} H. Cartan, {\it Vari\'et\'es analytiques r\'eelles et
  vari\'et\'es analytiques complexes}
  Bull. Soc. Math. France  {\bf 85}  (1957), 77--99.   

\bibitem{ca} H. Cartan, {\it Calculus diff\'erentiel.} Paris, Hermann
  1967. 

\bibitem{c} E. Chirka, {\sl Complex analytic sets.} Kluwer, Dordrecht,
  1989. 

\bibitem{df} K. Diederich and J. E. Forn\ae ss. {\it Pseudoconvex
  Domains with Real-Analytic Boundary.} Ann. Math. {\bf 107} (1978)
  371-384. 

\bibitem{dp} K. Diederich and S. Pinchuk, {\it The inverse of a
  CR-homeomorphism is CR}. Internat. J. Math. 4 (1993), no. 3,
  379--394. 

\bibitem{hi} H. Hironaka, {\it Subanalytic sets.}  Number theory,
  algebraic geometry and commutative algebra, in honor of Yasuo
  Akizuki,  pp. 453--493. Kinokuniya, Tokyo, 1973.   

\bibitem{l} S. Lojasiewicz, {\it Introduction to complex analytic
  geometry.} Birkh\"auser, Basel, 1991. 

\bibitem{m} B. Malgrange, {\it Ideals of differentiable functions.}
  Oxford Univ. Press, 1966.

\bibitem{n} R. Narasimhan, {\it Introduction to the Theory of Analytic
  Spaces.} Lecture Notes in Mathematics, vol. 25, Springer, 1966. 


\bibitem{ta} M. Tamm {\it Subanalytic sets in the calculus of
  variation.}  Acta Math. {\bf 146}  (1981), no. 3-4, 167--199. 

\bibitem{t} G. Tomassini,{\it Tracce delle funzioni olomorfe sulle
  sottovariet\'a analitiche reali d'una variet\'a complessa.}
  Ann. Scuola Norm. Sup. Pisa {\bf 20} (1966), 31--43.

\end{thebibliography}
}
\end{document}